\documentclass[11pt,a4paper]{article}

\usepackage{makeidx}
\makeindex            

\usepackage{amsmath}  
\usepackage{amssymb}  
\usepackage{amsthm}   
\usepackage{amsfonts}
\newcommand{\N}{\mathbb{N}}

\usepackage{mathrsfs}  

\newcommand{\Hom}[0]{\operatorname{Hom}}

\newcommand{\soc}[0]{\operatorname{soc}}
\newcommand{\Supp}[0]{\operatorname{Supp}}
\newcommand{\Ass}[0]{\operatorname{Ass}}

\newcommand{\height}[0]{\operatorname{height}}

\newtheorem{satz}{Theorem}[section]

\newtheorem{example}[satz]{Example}

\newtheorem{remark}[satz]{Remark}

\newtheorem{properties}[satz]{Properties}

\author{M. Hellus and J. St\"uckrad}
\title{Generalization of an example of Hartshorne
concerning local cohomology}
\date{\today}   

\begin{document}

\maketitle

\begin{abstract}
We prove the following generalization of an example of Hartshorne:
Let $k$ be a field, $n\geq 4$, $R=k[[X_1,\dots ,X_n]]$,
$I=(X_1,\dots ,X_{n-2})R$ and $p\in R$ a prime element such that
$p\in (X_{n-1},X_n)R$. Then $H^{n-2}_I(R/pR)$ is not artinian.
\end{abstract}
\section{Introduction}
It is an interesting question to determine if a given local
cohomology module $H^i_I(M)$ is artinian, where $I$ is an ideal of a
noetherian local ring $(R,m)$ and $M$ is a finite $R$-module; in
fact, this is one of Huneke's problems on local cohomology (see
\cite[third problem]{huneke92}).
\par
In general a given $R$-module $A$ is artinian if and only if both
\[ \Supp _R(A)\subseteq \{ m\} \tag{$A_1$}\]
and
\[ \soc (A):=\Hom _R(R/m,A)\hbox { is finitely generated}\tag{$A_2$}\] hold. Hartshorne (see \cite[section 3]{hartshorne70}) showed that
$H^2_I(R)$ is not artinian, where $k$ is a field,
$R=k[[u,v,x,y]]/(ux+vy)$ and $I\subseteq R$ is the ideal generated
by the classes of $u$ and $v$ in $R$ (in fact Hartshorne showed
something slightly different, but it is not difficult to modify his
result to the statement mentioned above). More precisely, $H^2_I(R)$
fulfills $(A_1)$, but not $(A_2)$.
\par
On the other hand, recently, modules of the form $D(H^i_I(M))$ have
been studied for general $R$, $I$, $i$ and $M$ (see
\cite{hellus05,hellusmArithm,hellusAttached,hellusStruct,hellusMatlisTop},
$D(\_ )$ is the Matlis dual functor, see \cite{matlis58} for details
on Matlis duality). It turns out that the study of these objects
leads to an elegant proof of Hartshorne's example; in fact we show
that $D(H^2_I(R))$ is not noetherian, where $H^2_I(R)$ is the local
cohomology module of Hartshorne's example. Hence $H^2_I(R)$ cannot
be artinian. Furthermore, our method can be used to show the
following more general result:
\begin{satz}
Let $k$ be a field, $R=k[[X_1,\dots ,X_n]]$ a power series algebra
over $k$ ($n\geq 3$), $I=(X_1,\dots ,X_{n-2})R$ and $p$ a prime
element of $R$ such that $p\in (X_{n-1},X_n)R$. Then
$H^{n-2}_I(R/pR)$ fulfills $(A_1)$, but not $(A_2)$; in particular,
it is not artinian.
\end{satz}It should be remarked that Marley and Vassilev (see \cite[theorem 2.3]{marley04}) have generalized Hartshorne's example in a different
direction. Due to different hypothesis, their and our generalization
can be compared only in a special case, see remark \ref{comparism}
for details.
\par
The authors thank Gennady Lyubeznik for drawing their attention to
Hartshorne's example.
\parindent=0pt
\bigskip
\bigskip
\section{Results}
\bigskip
Let $k$ by a field, $n\geq 3$, $R=k[[X_1,\dots ,X_n]]$ a formal
power series algebra over $k$, $I$ the ideal $(X_1,\dots ,X_{n-2})R$
of $R$ and define \[ D:=D(H^{n-2}_I(R))\ \ ,\] where $D(\_ )$ is the
Matlis dual functor, i. e. $D(M):=\Hom _R(M,E_R(R/m))$ for any
$R$-module $M$, $E_R(R/m)$ being a fixed $R$-injective hull of
$R/m$. Our method is to study the module $D$, here are some
properties of $D$:
\begin{properties}
(i) Every associated prime $p\in \Ass _R(D)$ has $\height (p)\leq
2$.
\par
(ii) For every prime ideal $p$ of $R$ with $\height (p)=2$ one has:
\[ p\in \Ass _R(D)\iff I+p \hbox { is }m\hbox
{-primary.}\]
\end{properties}
{\it Proof.} (i) Let $p\in \Ass _R(D)$. We conclude
\begin{eqnarray*}
0&\neq &\Hom _R(R/p,D)\\
&=&D(H^{n-2}_I(R)/pH^{n-2}_I(R))\\
&=&D(H^{n-2}_I(R/p))\ \ .\\
\end{eqnarray*}
Here the first equality follows formally from the exactness of $D(\_
)$ and the second from the right exactness of the functor
$H^{n-2}_I(\_ )=H^{n-2}_{(X_1,\dots ,X_{n-2})R}(\_ )$. It is
well-known that $H^{n-2}_I(R/p)\neq 0$ implies $\dim (R/p)\geq n-2$.
\par
(ii) Let $p$ be a prime ideal of $R$ of height two. Because of (i),
we have \[ p\in \Ass _R(D)\iff \Hom _R(R/p,D)\neq 0\ \ .\] In the
proof of (i) it was shown that the latter module equals \[
D(H^{n-2}_I(R/p)\] and, by the well-known Hartshorne-Lichtenbaum
vanishing theorem, is non-zero if and only if the ideal $I+p$ is
$m$-primary.
\begin{example}
\label{exampleH}In the above situation, take $n=4$ and, for every
$\lambda \in k$, define \[ p_\lambda :=(X_3+\lambda X_1,X_4+\lambda
X_2)R\ \ .\] Clearly, every $p_\lambda $ is a height two prime ideal
of $R$ and, by property (ii) above, is associated to
$D=D(H^2_I(R))$. On the other hand, for every $\lambda \in k$, one
has
\[ p:=X_1X_4+X_2X_3\in p_\lambda \] (because of
$p=X_1(X_4-\lambda X_2)+X_2(X_3+\lambda X_1)$). Therefore, at least
if $k$ is infinite, $D$ has infinitely many associated primes
containing $p$. This implies that \[ \Hom _R(R/pR,D)\] cannot be
finitely generated. But, as we have seen in the proof of property
(i) above, $\Hom _R(R/pR,D)$ is the Matlis dual of \[ H^2_I(R/pR)\]
and so $H^2_I(R/pR)$ cannot be artinian.
\end{example}
\begin{remark}
{\sloppy This is essentially Hartshorne's example (\cite[section
3]{hartshorne70}), the main difference is that Hartshorne works over
the ring $k[X_3,X_4][[X_1,X_2]]$, while we work over the ring
$k[[X_1,X_2,X_3,X_4]]$; but the two versions are essentially the
same, because the module
\[ H^2_{(X_1,X_2)}(k[X_3,X_4][[X_1,X_2]]/(X_1X_4+X_2X_3))\] is
naturally a module over $k[[X_1,X_2,X_3,X_4]]$, because its support
is $\{ (X_1,X_2,X_3,X_4)\} $. This is true, because for every prime
ideal $p$ of $k[X_3,X_4][[X_1,X_2]]$ different from
$(X_1,X_2,X_3,X_4)$ and containing $X_1X_4+X_2X_3$ the ring
$(k[X_3,X_4][[X_1,X_2]]/(X_1X_4+X_2X_3))_p$ is regular, and so
Hartshorne-Lichtenbaum vanishing shows that
\[ H^2_{(X_1,X_2)}(k[X_3,X_4][[X_1,X_2]]/(X_1X_4+X_2X_3))_p=0\ .\]}
\end{remark}
\bigskip A similar technique like
in example \ref{exampleH} works to show that $H^{n-2}_I(R/pR)$ is
not artinian for $R=k[[X_1,\dots ,X_n]]$, $n\geq 4$, $p\in
(X_{n-1},X_n)R$ a prime element and $k$ an arbitrary field (may be
finite):
\begin{satz}
{\sloppy \label{notArtinian} Let $k$ be a field, $n\geq 4$,
$R=k[[X_1,\dots ,X_n]]$, $I=(X_1,\dots ,X_{n-2})R$ and $p\in R$ a
prime element such that $p\in (X_{n-1},X_n)R$. Set
$D:=D(H^{n-2}_I(R))$.
\par
(i) If $p\in R$ is a prime element such that $p\in
(X_{n-1},X_n)R\cap I$ holds, the set \[ \{ p\in \Ass _R(D)\vert p\in
p\}
\] is infinite.
\par (ii) If $p\in R$ is a prime element such that $p\in
(X_{n-1},X_n)R$, $H^{n-2}_I(R/pR)$ is not artinian.}
\end{satz}
{\it Proof.}{\sloppy (i) It is easy to see that there exist $f,g\in
I, f\not\in (X_{n-1},X_n)R$ and $l\geq 1$ such that
\[ p=X_n^lf+X_{n-1}g\] holds (note that $f$ is not zero because $p$ is prime). Let $m\in \N ^+$ be arbitrary.
We have \[ p=(X_n^l+X_1^mg)f+(X_{n-1}-X_1^mf)g\] and so $p\in
I_m:=(X_n^l+X_1^mg,X_{n-1}-X_1^mf)R$. The elements $X_1,\dots
,X_{n-2},X_n^l+X_1^mg,X_{n-1}-X_1^mf$ form a system of parameters of
$R$ and so, by properties (i) and (ii) from above, there exists a
$p_m\in \Ass _R(D(H^{n-2}_I(R)))$ containing $I_m$. $p_m$
necessarily has height two. For $m,m^\prime \in \N ^+,m\neq m^\prime
$ \[ \sqrt {I_m+I_{m^\prime }}=(X_1,X_n,X_{n-1})R\cap \sqrt
{(X_{n-1},X_n,f,g)R}\] holds; in particular, all primes containing
$I_m+I_{m^\prime }$ have height at least three. The statement
follows now from property (i).
\par
(ii) If $p\not\in I$, it is easy to see that \[ \Supp
_R(H^{n-2}_I(R/pR))={\cal V}(I+pR)\ \ ,\] the set of prime ideals of
$R$ containing $I+pR$, and so $H^{n-2}_I(R/pR)$ does not satisfy
$(A_1)$. We assume $p\in I$: If $H^{n-2}_I(R/pR)$ was artinian, its
dual $D(H^{n-2}_I(R/pR))$ would be finitely generated; but we have
seen before that, because of the exactness of $D$ and the
right-exactness of $H^{n-2}_I(\_ )$, $D(H^{n-2}_I(R/pR))=\Hom
_R(R/pR,D)$, and from (i) we know that the latter module is not
finitely generated.}
\begin{remark}
Marley and Vassilev have shown \vskip 0.3cm {\bf Theorem}
(\cite[theorem 2.3]{marley04})
\par
Let $(T,m)$ be a noetherian local ring of dimension at least two.
Let $R=T[x_1,\dots ,x_n]$ be a polynomial ring in $n$ variables over
$T$, $I=(x_1,\dots ,x_n)$, and $f\in R$ a homogenous polynomial
whose coefficients form a system of parameters for $T$. Then the
*socle of $H^n_I(R/fR)$ is infinite dimensional.
\vskip 0.3cm In their paper \cite{marley04}, Marley and Vassilev say
(in section 1) that Hartshorne's example is obtained by letting
$T=k[[u,v]],n=2$ and $f=ux+vy$; there is a slight difference between
the two situations that comes from the fact that Hartshorne works
over a ring of the form $k[x,y][[u,v]]$ while Marley and Vassilev
work over a ring of the form $k[[u,v]][x,y]$. The two rings are not
the same. But, as
\[ \Supp _R(H^2_{(u,v)}(R/(uy+vx)))=\{ (x,y,u,v)\} \] (both for
$R=k[x,y][[u,v]]$ and for $R=k[[u,v]][x,y]$), the local cohomology
module in question is (in both cases) naturally a module over
$k[[x,y,u,v]]$ and, therefore, both versions are equivalent, i. e.
the result of Marley and Vassilev is a generalization of
Hartshorne's example.
\end{remark}
\begin{remark}
\label{comparism}\cite[theorem 2.3]{marley04} and our theorem
\ref{notArtinian} are both generalizations of Hartshorne's example,
but, due to different hypotheses, they can only be compared in the
following special case: $k$ a field, $n\geq 4$, \[
R_0=k[[X_{n-1},X_n]][X_1,\dots ,X_{n-2}]\ \ ,\]
\[ R=k[[X_1,\dots ,X_n]]\ \ ,\] $I=(X_1,\dots ,X_{n-2})R$, $p\in R_0$
a homogenous element such that $p$ is prime as an element of $R$.
Then \cite[theorem 2.3]{marley04} says (implicitly) that \[
H^{n-2}_I(R/pR)\] is not artinian, if the coefficients of $p\in R_0$
in $k[[X_{n-1},X_n]]$ form a system of parameters in
$k[[X_{n-1},X_n]]$, while theorem \ref{notArtinian} says that the
same module is not artinian if none of these coefficients of $p$ is
a unit in $k[[X_{n-1},X_n]]$.
\end{remark}

{\sc Universit\"at Leipzig, Fakult\"at f\"ur Mathematik und
Informatik, Ma\-the\-ma\-tisches Institut, Augustusplatz 10/11,
D-04109 Leipzig}

{\it E-mail}:

hellus@math.uni-leipzig.de

stueckrad@math.uni-leipzig.de
\end{document}